\documentclass[11pt,a4paper,leqno,twoside]{amsart}
\usepackage{latexsym,amssymb,amsmath}
\input xy
\xyoption{all}

\def\CC{\mathbb C}

\def\RR{\mathbb R}
\def\HH{\mathbb H}
\def\AA{{\mathbb A}}

\def\OO{\mathbb O}

\def\11{\mathbf 1}
\def\PP{\mathbb P}

\def\e1{\varepsilon_1}
\def\e2{\varepsilon_2}
\def\e3{\varepsilon_3}

\def\P2{{\PP}^2}

\def\00{\underline{0}}
\def\J0{{\cal J}_3(\underline{0})}

\def\PJ0{\PP({\cal J}_3(\underline{0}))}

\def\e{\varepsilon}

\def\AP2{{\AA\PP}^2}
\def\RP2{{\RR\PP}^2}
\def\CP2{{\CC\PP}^2}
\def\HP2{{\HH\PP}^2}
\def\OP2{{\OO\PP}^2}

\newtheorem{theo}{Theorem}[section]
\newtheorem{coro}[theo]{Corollary}
\newtheorem{lemm}[theo]{Lemma}

\theoremstyle{remark}
\newtheorem{rema}[theo]{Remark}

\begin{document}
\title[Locally finite-dimensional central division
algebras]{Locally finite-dimensional central division algebras
over function fields of curves over m-local fields, criterion
for normality}


\author{Ivan D. Chipchakov}
\address{Institute of Mathematics and Informatics\\Bulgarian
Academy of Sciences\\Acad. G. Bonchev Str., Bl. 8\\1113 Sofia,
Bulgaria:\\ E-mail: chipchak@math.bas.bg}

\begin{abstract}
Let $K _{m}$ be an $m$-local field with an $m$-th residue field $K
_{0}$, for some integer $m > 0$, and let $K/K _{m}$ be a field
extension of transcendence degree $1$. The paper under review shows
that if $K _{0}$ is a field of finite Diophantine dimension ddim$(K
_{0})$, and $R$ is an associative locally finite-dimensional central
division $K$-algebra, then $R$ is a normally locally finite algebra
over $K$, that is, every nonempty finite subset $Y$ of $R$ is
contained in a finite-dimensional central $K$-subalgebra
$\mathcal{R}_{Y}$ of $R$.
\end{abstract}

\maketitle
\par
\medskip
\section{\bf Introduction and statement of the main result}

\medskip
Let $F$ be a field, $A$ an (associative) unital $F$-algebra, $A
^{\ast }$ its multiplicative group, and $Z(A)$ the centre of $A$.
By a subalgebra of $A$, we mean a subalgebra containing the unit
of $A$. We say that $A$ is a locally finite-dimensional (abbr.,
LFD) algebra over $F$ if its finite subsets generate
finite-dimensional $F$-subalgebras; $A$ is said to be a central
$F$-algebra if $Z(A) = F$; it is called normally locally finite
(abbr., NLF) if every finite subset $Y \subset A$ is included in a
finite-dimensional central $F$-subalgebra $\Theta _{A}$ of $A$.
Clearly, NLF-algebras over $F$ are LFD and central $F$-algebras.
The study of infinite-dimensional division NLF-algebras has been
initiated by K\"{o}the, who has pointed out that the question of
whether central division LFD-algebras over $F$ are NLF is open
(see \cite[page~27]{Ko}). As shown in \cite{Ch7}, the answer to
this question is, generally, negative. Specifically, it has been
deduced from \cite[Theorems~2.1, 3.1]{Ch7}, the Lang-Nagata-Tsen
theorem (see \cite{Na}), and the main result of \cite{Mat} that
the class of central division LFD-algebras over a purely
transcendental extension $K$ of an algebraically closed field $K
_{0}$ consists of NLF-algebras over $K$ if and only if the
transcendence degree trd$(K/K _{0})$ is finite. The affirmative
answer to the considered question in case trd$(K/K _{0})$ is
finite is due to the fact that then $K$ is a field of type $C
_{t}$, in the sense of Lang, where $t = {\rm trd}(K/K _{0})$, and
by \cite[Theorem~2.3]{Ch7}, central division LFD-algebras over a
field $F _{m}$ of type $C _{m}$, for an integer $m \ge 0$, are
necessarily NLF over $F _{m}$. Note also that, by
\cite[Theorem~2.3]{Ch7}, for each $m \in \mathbb{N}$, central
division LFD-algebras over $m$-local fields with finite or
virtually perfect and PAC (pseudo-algebraically closed) $m$-th
residue fields are NLF as well.
\par
\smallskip
The notion of an $m$-local field (in the sense of \cite{HHKr} and
\cite{Li11}), where $m \in \mathbb{N}$, can be defined as follows:
by a $1$-local field, we mean a complete discrete valued field,
and when $m \ge 2$, an $m$-th local field with an $m$-th residue
field $K _{0}$ is a complete field $K _{m}$ with respect to a
discrete valuation $w _{0}$, such that the residue field $K
_{m-1}$ of $(K _{m}, w _{0})$ is an $(m-1)$-th local field with an
$(m-1)$-th residue field $K _{0}$. We say that a field $F$ is
virtually perfect if char$(F) = q$, and in case $q
> 0$, the degree $[F\colon F ^{q}]$ is finite, where $F ^{q} =
\{\lambda ^{q}\colon \lambda \in F\}$. This holds, for example, if
$F$ is a field of type $C _{m}$ (or a $C _{m}$-field), for an
integer $m \ge 0$, that is, if every $F$-form $f$ (a homogeneous
polynomial $f \neq 0$ with coefficients in $F$) of degree deg$(f)$
in more than deg$(f) ^{m}$ variables has a nontrivial zero over
$F$. It is well-known that $F$ is a $C _{0}$-field if and only if
it is algebraically closed; also, finite fields have type $C
_{1}$, by the Chevalley-Warning theorem, and PAC perfect fields
have type $C _{2}$ (cf. \cite[Proposition~21.2.4]{FJ}, and
\cite[Theorem~21.3.6 and Remark~21.3.7]{FJ}, respectively). It is
easily verified that if $F$ is a $C _{m}$-field and char$(F) = q >
0$, then $[F\colon F ^{q}] \le q ^{m}$ (the $F$-form $\sum _{i=1}
^{q^{m'}} a _{i}X _{i} ^{q}$ has no nontrivial zero over $F$ if
$[F\colon F ^{q}] = q ^{m'}$ and $a _{i} \in F\colon i = 1, \dots
, q ^{m'}$, are linearly independent over $F ^{q}$). By definition,
the Diophantine dimension ddim$(F)$ of $F$ is finite and equal to $r$
if $r$ is the least integer $\ge 0$, for which $F$ is a $C
_{r}$-field; ddim$(F)$ is infinity if $F$ is not a $C _{r'}$-field,
for any $r' \in \mathbb{N}$. In view of the Lang-Nagata-Tsen theorem,
the class of fields of finite Diophantine dimensions is closed under
the formation of field extensions of finite transcendence degrees,
and by Greenberg's theorem (see \cite{Gr}), it is closed under taking
$\mu $-fold iterated Laurent (formal power) series fields, for each
$\mu \in \mathbb{N}$. Therefore, the considered class contains any
field $\mathbb{F}$ of finite Diophantine dimension together with its
finitely-generated extensions. This, applied to the case where
$\mathbb{F}$ is finite, indicates that each finitely-generated field
of nonzero characteristic has a finite Diophantine dimension.
\par
\medskip
The purpose of this paper is to generalize Theorem~2.3 (c)
of \cite{Ch7}, by proving the normality of central division
LFD-algebras in the following situation:
\par
\medskip
\begin{theo}
\label{theo1.1} Let $K _{m}$ an $m$-local field with an $m$-th
residue field $K _{0}$, for some $m \in \mathbb{N}$, and let {\rm
ddim}$(K _{0}) < \infty $ and $K/K _{m}$ be a field extension with
{\rm trd}$(K/K _{m}) = 1$. Then central division {\rm
LFD}-algebras over $K$ are {\rm NLF}.
\end{theo}
\par
\medskip
The proof of Theorem \ref{theo1.1} is obtained as a consequence of
the main result of the following section.
\par
\medskip
\section{\bf $\Phi _{\rm Br}$-fields and normality of their central
division LFD-algebras}
\par
\smallskip
Let $F$ be a field, $\mathbb{P}$ the set of prime numbers,
$\mathbb{P}_{F} = \{p \in \mathbb{P}\colon p \neq {\rm
char}(F)\}$, Br$(F)$ the Brauer group of $F$, $s(F)$ the class of
associative finite-dimensional central simple $F$-algebras, and
$d(F)$ the subclass of the division algebras lying in $s(F)$. It
is known that Br$(F)$ is an abelian torsion group (cf.
\cite[Sects. 12.5 and 14.4]{P}), whence, it decomposes into the
direct sum $\oplus _{p \in \mathbb{P}} {\rm Br}(F) _{p}$ of its
$p$-components Br$(F) _{p}$. For each $\nabla \in s(F)$, denote by
ind$(\nabla )$ the Schur index of $\nabla $, i.e. the degree
deg$(D _{\nabla })$ of the underlying (central) division
$F$-algebra $D _{\nabla }$ of $\nabla $, determined by
Wedderburn's structure theorem (see \cite[Sect. 3.5]{P}); also,
let exp$(\nabla )$ be the exponent of $\nabla $, i.e. the order of
its (Brauer) equivalence class $[\nabla ]$ as an element of
Br$(F)$. As shown by Brauer (cf. \cite[Sect. 14.4]{P}),
exp$(\nabla )$ divides ind$(\nabla )$, and the sets of prime
divisors of ind$(\Delta )$ and exp$(\nabla )$ coincide. The proof
of these relations allows to obtain the above-noted properties of
Br$(F)$ and to deduce Brauer's primary tensor product
decomposition theorem for every $D \in d(F)$ (cf. \cite[Sects.
13.4, 14.4]{P}). The description of index-exponent pairs over $F$
depends on the Brauer $p$-dimensions Brd$_{p}(F)$, $p \in
\mathbb{P}$ (in the sense of Auel, Brussel, Garibaldi and Vishne,
see \cite{ABGV}), defined for each $p \in \mathbb{P}$, as follows:
Brd$_{p}(F) = n(p)$ if $n(p)$ is the least integer $\ge 0$ for
which ind$(A _{p}) \mid {\rm exp}(A _{p}) ^{n(p)}$ whenever $A
_{p} \in s(F)$ and $[A _{p}] \in {\rm Br}(F) _{p}$; when such
$n(p)$ does not exist, we put Brd$_{p}(F) = \infty $.
\par
The absolute Brauer $p$-dimension abrd$_{p}(E)$\footnote{The
Brauer $p$-dimension, in the sense of \cite{PS}, means the same as
the absolute Brauer $p$-dimension in the present paper.} of $E$ is
defined to be the supremum of Brd$_{p}(E ^{\prime })$, taken over
the set of finite extensions $E ^{\prime }$ of $E$ in $E _{\rm
sep}$, for every $p \in \mathbb{P}$. Denote by $\Phi _{\rm Br}$
the class of those fields $E$ for which there are integers $m
_{p}(E)\colon p \in \mathbb{P}$, with $m _{p}(E) \ge {\rm
Brd}_{p}(E ^{\prime })$, for every finite extension $E ^{\prime
}/E$ and any $p \in \mathbb{P}$; the fields from this class are
called Brauer finite-dimensional (or, briefly, $\Phi _{\rm
Br}$-fields). Clearly, if $E$ is a $\Phi _{\rm Br}$-field, then
abrd$_{p}(E) \le m _{p}(E)$, for every $p \in \mathbb{P}$. The
following lemma shows that
\par\noindent
abrd$_{p}(E) = m _{p}(E)$ whenever $p \in \mathbb{P}_{E}$; also,
it ensures that $E$ is a $\Phi _{\rm Br}$-field, provided that it
is virtually perfect with abrd$_{p}(E)$ finite, for every $p \in
\mathbb{P}_{E}$. For a proof of this lemma, we refer the reader to
\cite[Lemma~3.7]{Ch7}.
\par
\medskip
\begin{lemm}
\label{lemm2.1} Let $E$ be a field of characteristic $q > 0$,
$\overline E$ an algebraic closure of $E$, and $E _{\rm ins}$ the
maximal purely inseparable extension of $E$ in $\overline E$.
Then:
\par
{\rm (a)} {\rm abrd}$_{p}(E) \ge {\rm Brd}_{p}(E ^{\prime })$, for
each $p \in \mathbb{P}_{E}$, and every finite extension $E
^{\prime }$ of $E$ in $\overline E$; in addition, {\rm
abrd}$_{p}(E _{\rm ins}) = {\rm abrd}_{p}(E)$;
\par
{\rm (b)} If $E$ is virtually perfect and $[E\colon E ^{q}] = q
^{\delta }$, then {\rm Brd}$_{q}(E ^{\prime }) \le \delta $, for
every finite extension $E ^{\prime }/E$;
\par
{\rm (c)} $E _{\rm ins}$ is a $\Phi _{\rm Br}$-field, provided
that so is $E$; the converse holds if $E$ is virtually perfect.
\end{lemm}
\par
\medskip
The main result of this section can be stated as follows:
\par
\medskip
\begin{theo}
\label{theo2.2} Let $K _{m}$ be an $m$-local field with an $m$-th
residue field $K _{0}$ of characteristic $q$, for some $m \in
\mathbb{N}$, and let $\mathcal{K}_{m}/K _{m}$ and
$\mathcal{K}_{0}/K _{0}$ be purely transcendental extensions with
{\rm trd}$(\mathcal{K}_{i}/K _{i}) = 1$, where $i = 0$ and $i =
m$. Suppose that $\mathcal{K}_{0}$ is a $\Phi _{\rm Br}$-field,
denote by $\mathbb{P}_{q}$ the set $\mathbb{P} \setminus \{q\}$,
and for each
\par\vskip0.04truecm\noindent
$p \in \mathbb{P}$, put $b _{p}(K _{0}) = {\rm max}\{{\rm
abrd}_{p}(\mathcal{K}_{0}), 1 + {\rm abrd}_{p}(K _{0})\}$ in case
abrd$_{p}(\mathcal{K}_{0})
> 0$,
\par\vskip0.04truecm\noindent
and $b _{p}(K _{0}) = 0$ if abrd$_{p}(\mathcal{K}_{0}) = 0$. Then
$\mathcal{K}_{m}$ is a virtually perfect $\Phi _{\rm Br}$-field,
and {\rm abrd}$_{p}(\mathcal{K}_{m}) \le b _{p}(K _{0}) + m$, for
every $p \in \mathbb{P}_{q}$; hence, every central division {\rm
LFD}-algebra over $\mathcal{K}_{m}$ is {\rm NLF}.
\end{theo}
\par
\medskip
The latter assertion of Theorem \ref{theo2.2} follows from the
former one, since central division {\rm LFD}-algebras over any
$\Phi _{\rm Br}$-field are NLF (see \cite[Theorem~3.1]{Ch7}). The
former assertion of Theorem \ref{theo2.2} is proved in Section 4.
Here we show that Theorem \ref{theo1.1} can be deduced from
Theorem \ref{theo2.2}. The class of $C _{m}$-fields is closed
under the formation of algebraic extensions, for each $m$, and by
Matzri's theorem (cf. \cite{Mat}), Brd$_{p}(E _{m}) \le p ^{m-1} -
1$, for every $C _{m}$-field $E$ and each $p \in \mathbb{P} _{q}$.
As noted in Section 1, $C _{m}$-fields are virtually perfect,
whence, by Matzri's inequalities and Lemma \ref{lemm2.1}, they are
$\Phi _{\rm Br}$-fields. Assuming now that $K$, $K _{m}$ and $K
_{0}$ are fields satisfying the conditions of Theorem
\ref{theo1.1}, and applying consecutively the Lang-Nagata-Tsen
theorem to $K _{0}$, and the former part of Theorem \ref{theo2.2}
to $K _{m}$, one concludes that $K$ is a $\Phi _{\rm Br}$-field.
This, combined with \cite[Theorem~3.1]{Ch7}, proves Theorem
\ref{theo1.1}.
\par
\medskip
\begin{rema}
\label{rema2.3} A well-known conjecture (stated by M. Artin in the
case of $d =2$) predicts that if $K$ is a $C _{d}$-field, for some
$d > 0$, then Brd$_{p}(K) \le d - 1$, for all $p \in \mathbb{P}$.
A proof of this conjecture will reduce Theorem \ref{theo1.1} to a
consequence of \cite[Theorem~3]{PS}, and results of \cite{HHKr}
and \cite{Li11} (generalizing Saltman's theorem proved in
\cite{Salt}). Theorem \ref{theo2.2} is a modified version of these
results, which allows to deduce Theorem \ref{theo1.1} from
\cite[Theorem~3.1]{Ch7}, independently of the generalized Artin
conjecture.
\end{rema}
\par
\medskip
The basic notation, terminology and conventions kept in this paper
are standard and essentially the same as in \cite{TW}, \cite{L},
\cite{P} and \cite{S1}. For convenience of the reader, we present
in Section 4 a part of the corresponding information concerning
fields with Henselian valuations. Brauer groups and ordered
abelian groups are written additively, Galois groups are viewed as
profinite with respect to the Krull topology, and by a profinite
group homomorphism, we mean a continuous one. Throughout,
$\mathbb{Z}$ is the additive group of integers, and for any field
$E$, $E _{\rm sep}$ is a separable closure of $E$. Given a field
extension $E ^{\prime }/E$, we write I$(E ^{\prime }/E)$ for the
set of intermediate fields of $E ^{\prime }/E$. When $E ^{\prime
}/E$ is Galois, $\mathcal{G}(E ^{\prime }/E)$ denotes its Galois
group, and $\mathcal{G}_{E}$ stands for the absolute Galois group
of $E$, i.e. $\mathcal{G}_{E} = \mathcal{G}(E _{\rm sep}/E)$. If
$\mathcal{F}/F$ is a purely transcendental field extension with
trd$(\mathcal{F}/F) = 1$ and abrd$_{p}(\mathcal{F})$ finite, for
some $p \in \mathbb{P}$, we put $b _{p}(F) = {\rm max}\{{\rm
abrd}_{p}(\mathcal{F}), 1 + {\rm abrd}_{p}(F)\}$ if
abrd$_{p}(\mathcal{F}) > 0$, and $b _{p}(F) = 0$ in case
abrd$_{p}(\mathcal{F}) = 0$.
\par
\medskip
\section{\bf Two lemmas on absolute Brauer $p$-dimensions}
\par
\medskip
This section contains lemmas which play a role in our proof of
Theorem \ref{theo2.2}. The first lemma presents Galois-theoretic
ingredients of our proof.
\par
\medskip
\begin{lemm}
\label{lemm3.1} Let $E$ be a field, $p$ a prime number, $G _{p}$ a
Sylow pro-$p$ subgroup of $\mathcal{G}_{E}$, and $E _{p}$ the
fixed field of $G _{p}$; also, let $\mathcal{E}$ be an extension
of $E$ in $E _{\rm sep}$, such that $p \nmid [\widetilde{E}
_{0}\colon E]$, for any finite extension $\widetilde{E}_{0}$ of
$E$ in $\mathcal{E}$. Then {\rm Brd}$_{p}(E) \le {\rm
Brd}_{p}(\mathcal{E})$ and {\rm abrd}$_{p}(E) = {\rm
abrd}_{p}(\mathcal{E}) = {\rm abrd}_{p}(E _{p})$.
\end{lemm}
\par
\medskip
\begin{proof}
For any field $E ^{\prime } \in I(E _{\rm sep}/E)$, denote by
$\Lambda _{E'}$ the set of those
\par\noindent
$E _{1} ^{\prime } \in I(E _{\rm sep}/E ^{\prime })$, for which
$p$ does not divide the degree of any finite extension of $E
^{\prime }$ in $E _{1} ^{\prime }$. Clearly, $E ^{\prime } \in
\Lambda _{E'}$, i.e. $\Lambda _{E'} \neq \emptyset $, and viewed
as a set partially ordered by inclusion, $\Lambda _{E'}$ satisfies
the conditions of Zorn's lemma, whence, it contains a maximal
element. Note also that $\Lambda _{E''} = I(E _{\rm sep}/E
^{\prime \prime }) \cap \Lambda _{E'}$, for every $E ^{\prime
\prime } \in \Lambda _{E'}$. Since finite extensions of $E
^{\prime }$ in $E _{\rm sep}$ are simple (cf. \cite[Ch.~V,
Theorem~4.6]{L}), this follows from the fact that, for each $\beta
\in E _{\rm sep}$, there is a finite extension $E ^{\prime }
_{\beta }$ of $E ^{\prime }$ in $E ^{\prime \prime }$, such that
$[E ^{\prime } _{\beta }(\beta )\colon E ^{\prime } _{\beta }] =
[E ^{\prime \prime }(\beta )\colon E ^{\prime \prime }]$ (whence,
$p \mid [E _{\beta } ^{\prime }(\beta )\colon E ^{\prime }]$ if
and only if $p \mid [E _{\beta } ^{\prime }(\beta )\colon E
_{\beta } ^{\prime }]$); one may take as $E ^{\prime } _{\beta }$
the extension of $E ^{\prime }$ generated by the coefficients of
the minimal (monic) polynomial of $\beta $ over $E ^{\prime \prime
}$. These observations allow to deduce from Galois theory and
Sylow's theorems (more precisely, from their version for profinite
groups, see \cite[Ch. I, 1.4]{S1}) that a field $\mathcal{E}
^{\prime } \in I(E _{\rm sep}/E ^{\prime })$ is a maximal element
in $\Lambda _{E'}$ if and only if $\mathcal{G}_{\mathcal{E}'}$ is
a Sylow pro-$p$ subgroup of $\mathcal{G}_{E'}$. Since $\mathcal{E}
\in \Lambda _{E}$, it is now easy to see that a closed subgroup
$\mathcal{H} _{p}$ of $\mathcal{G}_{\mathcal{E}}$ is a Sylow
pro-$p$ subgroup of $\mathcal{G}_{E}$ if and only if it is a Sylow
pro-$p$ subgroup of $\mathcal{G}_{\mathcal{E}}$. Thus it turns out
that Lemma \ref{lemm3.1} will be proved if we show that
Brd$_{p}(E) \le {\rm Brd}_{p}(\mathcal{E})$ and abrd$_{p}(E) =
{\rm abrd}_{p}(E _{p})$. The former inequality follows from the
assumption on $\mathcal{E}$, which ensures that the scalar
extension map $s(E) \to s(\mathcal{E})$ induces an
index-preserving injective group homomorphism Br$(E) _{p} \cong
{\rm Br}(\mathcal{E}) _{p}$.
\par
We turn to the proof of the equality abrd$_{p}(E) = {\rm
abrd}_{p}(E _{p})$. It is easy to see that, for each pair $Y \in
I(E _{\rm sep}/E)$, $\Delta _{Y} \in d(Y)$, there exists a pair
\par\vskip0.034truecm\noindent
$Y _{0} \in I(Y/E)$, $\Delta _{Y _{0}} \in d(Y _{0})$, such that
$[Y _{0}\colon E] < \infty $, exp$(\Delta _{Y _{0}}) = {\rm
exp}(\Delta _{Y})$,
\par\vskip0.034truecm\noindent
and there is an $Y$-isomorphism $\Delta _{Y _{0}} \otimes _{Y
_{0}} Y \cong \Delta _{Y}$ (cf. \cite[(1.3)]{Ch2}); in particular,
deg$(\Delta _{Y _{0}}) = {\rm deg}(\Delta _{Y})$. This, applied to
the case where $Y$ is an arbitrary finite extension of $E _{p}$
and deg$(\Delta _{Y})$ is a $p$-power, leads to the conclusion
that Brd$_{p}(Y) \le {\rm abrd}_{p}(E)$ and abrd$_{p}(E _{p}) \le
{\rm abrd}_{p}(E)$.
\par
Our objective now is to prove that abrd$_{p}(E) \le {\rm
abrd}_{p}(E _{p})$. It is clearly sufficient to show that given a
finite extension $Z$ of $E$ in $E _{\rm sep}$, and an algebra $D
\in d(Z)$ of $p$-power degree, there is a finite extension $Z
_{p}$ of $E _{p}$ in $E _{\rm sep}$, such that exp$(D _{p}) = {\rm
exp}(D)$ and deg$(D _{p}) = {\rm deg}(D)$, for some $D _{p} \in
d(Z _{p})$. Denote by $M$ the Galois closure of $Z$ in $E _{\rm
sep}$ over $E$, fix a Sylow $p$-subgroup $H _{p}$ of
$\mathcal{G}(M/Z)$, choose a Sylow $p$-subgroup $G _{p}$ of
$\mathcal{G}(M/E)$ so that $H _{p} \le G _{p}$, and let $Z
^{\prime }$ and $U$ be the fixed fields of $H _{p}$ and $G _{p}$,
respectively. It follows from Galois theory and the choice of $H
_{p}$ and $G _{p}$ that $Z.U \subseteq Z ^{\prime }$, $p$ does not
divide $[Z ^{\prime }\colon Z][U\colon E]$, and $[Z ^{\prime
}\colon U]$ equals both the index of $H _{p}$ in $G _{p}$, and the
maximal $p$-power dividing $[Z\colon E]$. These calculations yield
consecutively $[Z ^{\prime }\colon E]$ and $Z ^{\prime } = Z.U$.
They also ensure that $d(Z ^{\prime })$ contains the $Z ^{\prime
}$-algebra $D ^{\prime } = D \otimes _{Z} Z ^{\prime }$, and prove
that exp$(D ^{\prime }) = {\rm exp}(D)$ and deg$(D ^{\prime }) =
{\rm deg}(D)$.
\par
Assume now that $\widetilde G _{p}$ is a Sylow pro-$p$ subgroup of
$\mathcal{G}_{U}$, $U _{p}$ is the fixed field of $\widetilde G
_{p}$, and $Z ^{\prime } _{p} = Z ^{\prime }.U _{p}$. As noted
above, $\widetilde G _{p}$ is a Sylow pro-$p$ subgroup of
$\mathcal{G}_{E}$, and by the version of Sylow's theorems for
profinite groups, the groups $G _{p}$ and $\widetilde G _{p}$ are
conjugate in $\mathcal{G}_{E}$; hence, by Galois theory, there
exists $\tau _{p} \in \mathcal{G}_{E}$, which induces an
isomorphism $U _{p} \cong E _{p}$ as $E$-algebras. Observing now
that each finite extension $\nabla ^{\prime }$ of $Z ^{\prime }$
in $Z ^{\prime } _{p}$ is a subfield of the compositum $Z ^{\prime
}\nabla $, for some finite extension $\nabla $ of $U$ in $U _{p}$,
and we have $[\nabla Z ^{\prime }\colon Z ^{\prime }] = [\nabla
\colon U]$, one obtains that $p \nmid [\nabla ^{\prime }\colon Z
^{\prime }]$. This implies the $Z _{p} ^{\prime }$-algebra $D _{p}
^{\prime } = D ^{\prime } \otimes _{Z'} Z _{p} ^{\prime }$ lies
\par\vskip0.031truecm\noindent
in $d(Z _{p} ^{\prime })$, and gives exp$(D _{p} ^{\prime }) =
{\rm exp}(D)$ and deg$(D _{p} ^{\prime }) = {\rm deg}(D)$. At the
same
\par\vskip0.031truecm\noindent
time, it follows that $Z _{p} ^{\prime }/U _{p}$ is a field
extension with $[Z _{p} ^{\prime }\colon U _{p}] = [Z ^{\prime
}\colon U]$, $Z _{p} ^{\prime }$ is the fixed field of the group
$\widetilde H _{p} = \widetilde G _{p} \cap \mathcal{G}_{Z'}$, and
$\widetilde H _{p}$ is a Sylow pro-$p$ subgroup of
$\mathcal{G}_{Z'}$. Finally, one concludes that if $Z _{p}$ is the
image of $Z _{p} ^{\prime }$ under $\tau _{p}$, then $Z _{p}/E
_{p}$ is a field extension, $[Z _{p}\colon E _{p}] = [Z _{p}
^{\prime }\colon U _{p}] = [Z ^{\prime }\colon U]$, and $\tau
_{p}$ gives rise
\par\vskip0.031truecm\noindent
to an index-preserving isomorphism $\bar \tau _{p}: {\rm Br}(Z
_{p} ^{\prime }) _{p} \cong {\rm Br}(Z _{p}) _{p}$. Therefore, if
\par\vskip0.031truecm\noindent $D _{p} \in d(Z _{p})$ satisfies
$[D _{p}] = \bar \tau _{p} ([D _{p} ^{\prime }])$, then exp$(D
_{p}) = {\rm exp}(D _{p} ^{\prime }) = {\rm exp}(D ^{\prime })$
\par\vskip0.031truecm\noindent
and deg$(D _{p}) = {\rm deg}(D _{p} ^{\prime }) = {\rm deg}(D
^{\prime })$. This completes the proof of the inequality
abrd$_{p}(E) \le {\rm abrd}_{p}(E _{p})$, so Lemma \ref{lemm3.1}
is proved.
\end{proof}
\par
\smallskip
\begin{lemm}
\label{lemm3.2} Assuming that $E$, $p$, $G _{p}$ and $E _{p}$
satisfy the conditions of Lemma \ref{lemm3.1}, put $\mathbb{P}
^{\prime } = \{p' \in \mathbb{P}_{E}\colon p' \neq p\}$, $F =
E(X)$ and $F _{1} = E _{\rm sep}(X)$, where $X$ is a
transcendental element over $E$. Then {\rm abrd}$_{p}(E) \le {\rm
Brd}_{p}(F)$, {\rm abrd}$_{p}(E _{p}(X)) = {\rm abrd}_{p}(F)$, and
{\rm abrd}$_{p'}(E _{p}(X)) = 0$, for every $p' \in \mathbb{P}
^{\prime }$.
\end{lemm}
\par
\smallskip
\begin{proof}
The inequality abrd$_{p}(E) \le {\rm Brd}_{p}(F)$ is a special
case of Theorem~2.1 (a) and (b) of \cite{Ch3}. The proof of the
equality abrd$_{p}(E _{p}(X)) = {\rm abrd}_{p}(F)$ relies on the
fact that $E$ has no proper algebraic extensions in $F$, which
allows to deduce from Galois theory (cf. \cite[Ch. VI,
Theorem~1.12]{L}) that there exists a canonical group isomorphism
$\varphi \colon \mathcal{G}_{E} \to \mathcal{G}(F _{1}/F)$, such
that the image of $\mathcal{G}_{E'}$ under $\varphi $ coincides
with $\mathcal{G}(F _{1}/E ^{\prime }(X))$, for each $E ^{\prime }
\in I(E _{\rm sep}/E)$. At the same time, it becomes clear that
the mapping of $I(E _{\rm sep}/E)$ into $I(F _{1}/E)$, by the rule
$E ^{\prime } \to E ^{\prime }(X)$, is a degree-preserving and
(set-theoretic) inclusion-preserving bijection. Therefore, it is
not difficult to see that $\mathcal{G}(F _{1}/E _{p}(X))$ is a
Sylow pro-$p$ subgroup of $\mathcal{G}(F _{1}/F)$; in particular,
$p$ does not divide the degree of any finite extension of $F$ in
$E _{p}(X)$, which reduces the equality abrd$_{p}(E _{p}(X)) =
{\rm abrd}_{p}(F)$ to a consequence of Lemma \ref{lemm3.1}. Now
fix an algebraic closure $\overline E$ of $E$, and an arbitrary
finite extension $E _{p} ^{\prime }$ of $E _{p}(X)$ in $F _{\rm
sep}$, and put $\widetilde E _{p} = E _{p} ^{\prime }F _{1}$. We
show that abrd$_{p'}(E _{p}(X)) = 0$, for all $p' \in \mathbb{P}
^{\prime }$, proving in conjunction with the Albert-Hochschild
theorem (see \cite[Ch.~II, 2.2]{S1}) that Br$(E _{p} ^{\prime })
_{p'} = \{0\}$. It is easily verified that $\widetilde E _{p}/F
_{1}$ is a finite extension, so it follows from Tsen's theorem
that Br$(\widetilde E _{p}) = \{0\}$ in case char$(E) = 0$. If
char$(E) = q > 0$, then $\widetilde E _{p}F _{1,{\rm ins}}/F
_{1,{\rm ins}}$ is a finite extension, $\widetilde E _{p}F
_{1,{\rm ins}} = \widetilde E _{p,{\rm ins}}$, and the set $I(F
_{1,{\rm ins}}/F _{1})$ contains an isomorphic copy of the field
$\overline E(X)$; in view of Tsen's theorem, this yields
Br$(\widetilde E _{p,{\rm ins}}) = \{0\}$. Hence, by Lemma
\ref{lemm2.1} (a), the assertion that Br$(\widetilde E _{p}) _{p'}
= \{0\}$ holds in general, for all $p' \in \mathbb{P}_{E}$, which
means that $\widetilde E _{p}$ is a splitting field of every
$\Delta \in d(E _{p} ^{\prime })$ of degree not divisible by
char$(E)$. Our conclusion can be restated by saying that $\Delta $
has a splitting field $E _{p} ^{\prime }(\Delta )$ which is a
finite extension of $E _{p} ^{\prime }$ in $\widetilde E _{p}$
(see \cite[Sect. 13.2, Proposition~b]{P} or \cite[(1.3)]{Ch2}).
Observing finally that $\widetilde E _{p}/E _{p} ^{\prime }$ is a
Galois extension with $\mathcal{G}(\widetilde E _{p}/E _{p}
^{\prime })$ a pro-$p$ group, one obtains that finite extensions
of $E _{p} ^{\prime }$ in $\widetilde E _{p}$ are of $p$-power
degrees. It is now easy to show that if char$(E) \nmid {\rm
deg}(\Delta )$, then deg$(\Delta )$ is a $p$-power, which implies
Br$(E _{p} ^{\prime }) _{p'} = \{0\}$, for all $p' \in \mathbb{P}
^{\prime }$, and so completes the proof of Lemma \ref{lemm3.2}.
\end{proof}
\par
\smallskip
\section{\bf Proof of Theorem \ref{theo1.1}}
\par
\smallskip For any field $K$ with a (nontrivial) Krull valuation
$v$, $v(K)$ denotes the value group of $(K, v)$, $O _{v}(K) = \{a
\in K\colon \ v(a) \ge 0\}$ is the valuation ring of $(K, v)$, $M
_{v}(K) = \{\mu \in K\colon \ v(\mu ) > 0\}$ the maximal ideal of
$O _{v}(K)$, and $\widehat K = O _{v}(K)/M _{v}(K)$ the residue
field of $(K, v)$. As usual, $v(K)$ is assumed to be an ordered
abelian group; also, $\overline {v(K)}$ stands for a fixed
divisible hull of $v(K)$. The valuation $v$ is said to be
Henselian if it extends uniquely, up-to equivalence, to a
valuation $v _{L}$ on each algebraic extension $L$ of $K$. When
this holds, $(K, v)$ is called a Henselian field. It is known that
$(K, v)$ is Henselian if it is maximally complete, i.e. it has no
valued extension $(K ^{\prime }, v')$ such that $K ^{\prime } \neq
K$, $\widehat K ^{\prime } = \widehat K$ and $v'(K ^{\prime }) =
v(K)$ (cf. \cite[Theorem~15.3.5]{E3}). \label{max_comp} For
example, complete discrete valued fields are maximally complete,
and for each $n \in \mathbb{N}$, so is the $n$-fold iterated
Laurent series field $K _{n} = K _{0}((X _{1})) \dots ((X _{n}))$
over a field $K _{0}$ (that is, $K _{m} = K _{m-1}((X _{m}))$, for
every index $m > 0$) with respect to its standard valuation $v
_{n}$ inducing on $K _{0}$ the trivial valuation (see
\cite[Exercise~3.11]{TW}, and \cite[Theorem~18.4.1]{E3}). Here $v
_{n}(K _{n}) = \mathbb{Z} ^{n}$, $K _{0}$ is the residue field of
$(K _{n}, v _{n})$, and the abelian group $\mathbb{Z} ^{n}$ is
considered with its inverse-lexicographic ordering (see
\cite[Examples~4.2.2 and 9.2.2]{E3}).
\par
When $v$ is Henselian, so is $v _{L}$, for any algebraic field
extension $L/K$. In this case, we denote by $\widehat L$ the
residue field of $(L, v _{L})$, and put $v(L) = v _{L}(L)$, $O
_{v}(L) = O _{v _{L}}(L)$. Clearly, $\widehat L/\widehat K$ is an
algebraic extension and $v(K)$ is an ordered subgroup of $v(L)$,
such that $v(L)/v(K)$ is a torsion group; hence, one may assume
without loss of generality that \label{k99} $v(L)$ is an ordered
subgroup of $\overline {v(K)}$. By Ostrowski's theorem (cf.
\cite[Theorem~17.2.1]{E3}), if $[L\colon K]$ is finite, then it is
divisible by $[\widehat L\colon \widehat K]e(L/K)$, where $e(L/K)$
is the index of $v(K)$ in $v(L)$. The equality $[L\colon K] =
[\widehat L\colon \widehat K]e(L/K)$ holds if char$(\widehat K)
\nmid [L\colon K]$ (apply Ostrowski's theorem), or else, if $v$ is
discrete and $L/K$ is separable (see \cite[Sect. 17.4]{E3}). We
say that $L/K$ is inertial if $[L\colon K] = [\widehat L\colon
\widehat K]$ and $\widehat L/\widehat K$ is separable (when this
holds, $L/K$ is separable). Inertial extensions have the following
properties (for a proof, see \cite[Theorem~A.23 and
Corollary~A.25]{TW}):
\par
\medskip
\begin{lemm}
\label{lemm4.3} Let $(K, v)$ be a Henselian field and $K _{\rm
ur}$ the compositum of inertial extensions of $K$ in $K _{\rm
sep}$. Then:
\par
{\rm (a)} An inertial extension $R ^{\prime }/K$ is Galois if and
only if so is $\widehat R ^{\prime }/\widehat K$. When this holds,
$\mathcal{G}(R ^{\prime }/K)$ and $\mathcal{G}(\widehat R ^{\prime
}/\widehat K)$ are canonically isomorphic.
\par
{\rm (b)} $v(K _{\rm ur}) = v(K)$, $K _{\rm ur}/K$ is a Galois
extension and $\mathcal{G}(K _{\rm ur}/K) \cong
\mathcal{G}_{\widehat K}$.
\par
{\rm (c)} Finite extensions of $K$ in $K _{\rm ur}$ are inertial, and
the natural mapping of $I(K _{\rm ur}/K)$ into $I(\widehat K _{\rm
sep}/\widehat K)$, by the rule $L \to \widehat L$, is bijective.
\par
{\rm (d)} For each finite extension $K _{1}$ of $K$ in $K _{\rm
sep}$, the field $K _{0} = K _{1} \cap K _{\rm ur}$ equals the
maximal inertial extension of $K$ in $K _{1}$; in addition,
$\widehat K _{0} = \widehat K _{1}$.
\end{lemm}
\par
\medskip
Our next lemma, combined with Lemmas \ref{lemm3.1} and
\ref{lemm3.2}, makes it possible to prove Theorem \ref{theo2.2},
by applying the Harbater-Hartmann-Krashen and Lieblich theorem
(see \cite[Corollaries~5.7, 5.8]{HHKr} and
\cite[Corollary~1.3]{Li11}). For a  presentation of excellent
rings (introduced in the theory of Noetherian commutative rings by
Grothendieck), we refer the reader to \cite{Li}.
\par
\medskip
\begin{lemm}
\label{lemm4.4} Let $E$ be a virtually perfect field with char$(E)
= q > 0$ and $[E\colon E ^{q}] = q ^{\delta }$, and let $U$ be an
extension of $E$ in $E _{\rm sep}$, $E _{1} = E((Z))$ and $U _{1}
= U((Z))$ be Laurent series fields in an indeterminate $Z$ over
$E$ and $U$, respectively. Put $W = U.E _{1}$, and denote by $e$
and $u$ the standard (Henselian) $\mathbb{Z}$-valued valuations of
$E _{1}$ and $U _{1}$, trivial on $E$ and $U$, respectively. Then
the ring $O _{w}(W)$, where $w = e _{W}$, is excellent.
\end{lemm}
\par
\smallskip
\begin{proof} Put $B _{1} = B \cup \{Z\}$, where $B$ is a $q$-basis
of $E$ over $E ^{q}$ or $B = \emptyset $ depending on whether or
not $E \neq E ^{q}$. It is known that $U = U ^{q}$ if $E = E
^{q}$, and $B$ is a $q$-basis of $U$ over $U ^{q}$ in case $E \neq
E ^{q}$ (cf. \cite[Lemma~2.7.3]{FJ}); in particular, $[U\colon U
^{q}] = q ^{\delta }$. Note also that $E _{1} ^{q} = E ^{q}((Z
^{q}))$, $U _{1} ^{q} = U ^{q}((Z ^{q}))$, and $W ^{q} = U ^{q}.E
_{1} ^{q}$, which implies $B _{1}$ is a $q$-basis of $E _{1}$, $U
_{1}$ and $W$ over $E _{1} ^{q}$, $U _{1} ^{q}$ and $W ^{q}$,
respectively. At the same time, it is easily verified that
\par\noindent
$O _{e}(E _{1}) = E[[Z]]$, $O _{u}(U _{1}) = U[[Z]]$, and $O
_{w}(W) = U.E[[Z]] = U ^{q}.E[[Z]]$, where $w = e _{W}$. This
indicates that $[W\colon W ^{q}] = q ^{\delta }$, $u$ is the
continuous prolongation of $w$ upon $U _{1}$, the set $O _{w}(W)
^{q} = \{r ^{q}\colon r \in O _{w}(W)\}$ is a subring of $O
_{w}(W)$ equal to $U ^{q}.E ^{q}[[Z ^{q}]]$, and $O _{w}(W)$ is a
free module over $O _{w}(W) ^{q}$ of rank $q ^{\delta +1}$.
Therefore, by the Kunz theorem \cite{Kunz}, $O _{w}(W)$ is an
excellent ring. Equivalently, one obtains that $U _{1}/W$ is a
separable extension, in the sense of \cite{L} and
\cite[Lemma~2.6.1]{FJ}, which is algebraic if and only if
$[U\colon F] < \infty $.
\end{proof}
\par
We are now prepared to complete the proof of Theorem
\ref{theo2.2}. Put
\par\noindent
$K = \mathcal{K}_{m} = K _{m}(X)$, where $X$ is an indeterminate
over $K _{m}$. We show that
\par\noindent
$b _{p}(K _{m}) \le b _{p}(K _{0}) + m$, for each $p \in
\mathbb{P}_{q}$, and in case $q > 0$, abrd$_{q}(\mathcal{K}_{m}) <
\infty $.
\par
Suppose first that $m = 1$, char$(K _{1}) = 0$, and char$(K _{0}) = q
> 0$. By definition, $K _{1}$ is complete with respect to a discrete
valuation $v _{1}$ such that $\widehat K _{1} = K _{0}$. Note
further that $\mathcal{K}_{0} \in \Phi _{\rm Br}$, which ensures
that $K _{0}$ is a virtually perfect field (cf.
\cite[Theorem~2.1~(c)]{Ch3}), so it follows from
\cite[Theorem~3]{PS} that abrd$_{q}(K) < \infty $. Henceforth, we
assume that $p \neq q$ and $K _{p}$ is the fixed field of some
Sylow pro-$p$ subgroup $G _{p}$ of $\mathcal{G}_{K}$. Consider the
fields $\Theta _{p} = K _{1,{\rm sep}}(X) \cap K _{p}$,
$\mathcal{U}_{p} = K _{1,{\rm ur}}(X) \cap K _{p}$, $K _{1,p} =
\Theta _{p} \cap K _{1,{\rm sep}}$ and $U _{p} = \mathcal{U} _{p}
\cap K _{1,{\rm ur}}$. Using Galois theory and the form of Sylow's
theorems for profinite groups as in the proof of Lemma
\ref{lemm3.2}, one obtains that $\Theta _{p} = K _{1,p}(X)$,
$\mathcal{U}_{p} = U _{p}(X)$, and $\mathcal{G}_{K _{1,p}}$ and
$\mathcal{G}(K _{1,{\rm ur}}/U _{p})$ are Sylow pro-$p$ subgroups
of $\mathcal{G}_{K _{1}}$ and $\mathcal{G}(K _{1,{\rm ur}}/K
_{1})$, respectively. Therefore, by Lemma \ref{lemm3.1},
abrd$_{p}(K) = {\rm abrd}_{p}(\Theta _{p}) = {\rm
abrd}_{p}(\mathcal{U}_{p})$, and by Lemma \ref{lemm4.3}, the
valuation of $U _{p}$ extending $v _{1}$ is Henselian and
discrete, $\widehat U _{p} \in I(K _{0,{\rm sep}}/K _{0})$, and
$\mathcal{G}_{\widehat U _{p}}$ is a Sylow pro-$p$ subgroup of
$\mathcal{G}_{K_{0}}$. Let now $\widehat U _{p}(X _{0})$ be a
rational function field in an indeterminate $X _{0}$ over
$\widehat U _{p}$. Then, by Lemma \ref{lemm3.2},
abrd$_{p}(\widehat U _{p}(X _{0})) = {\rm abrd}_{p}(K _{0}(X
_{0}))$ and abrd$_{p'}(\widehat U _{p}(X _{0})) = 0$, for every
$p' \in \mathbb{P}$ with $p \neq p' \neq q$. Note also that
abrd$_{p}(\widehat U _{p}) = {\rm abrd}_{p}(K _{0})$, by Lemma
\ref{lemm3.1}, which ensures that $b _{p}(\widehat U _{p}) = b
_{p}(K _{0})$. Taking into account that the ring $O _{v
_{1}}(\widehat U _{p})$ is excellent (see
\cite[Corollary~8.2.40]{Li}), one obtains from
\cite[Corollaries~5.7, 5.8]{HHKr} and \cite[Corollary~1.3]{Li11})
that $b _{p}(U _{p}) \le b _{p}(K _{0}) + 1$. Thus Theorem
\ref{theo2.2} is proved in the special case of $m = 1$ and char$(K
_{0}) = 0 \neq q$.
\par
Next it is obtained by induction on $m$ that, for the rest of our
proof, one may assume that char$(K _{m}) = {\rm char}(K _{0}) = q$
and $K _{0}$ is a virtually perfect field. In view of the
Hasse-Schmidt theorem (see the references on \cite[page~110]{E3}),
this means to consider only the case of $K _{m} = K _{0}((Y _{1}))
\dots ((Y _{m}))$. Then $K _{m-j} = K _{0}((Y _{1})) \dots ((Y
_{m-j}))$ is the $j$-th residue field of $K _{m}$, and it is
virtually perfect, for each index $j > 0$; hence, by Lemma
\ref{lemm2.1}, abrd$_{q}(K _{m})$ is finite when $q > 0$. As in
the case of $m = 1$, fix an arbitrary $p \in \mathbb{P}$, denote
by $K _{0,p}$ the fixed field of some Sylow pro-$p$ subgroup of
$\mathcal{G}_{K _{0}}$, and identifying $K _{0,{\rm sep}}$ with
its $K _{0}$-isomorphic copy in $K _{m,{\rm sep}}$, put $L _{m-j}
= K _{0,p}K _{m-j}$, $0 \le j \le m - 1$. Also, let $v _{j}$ be
the prolongation on $L _{m-j}$ of the standard $\mathbb{Z}$-valued
valuation of $K _{m-j}$ that is trivial on $K _{m-j-1}$, for each
index $j < m$. Clearly, $L _{m-j}$ is an algebraic extension of $K
_{m-j}$ generated by $K _{0,p}$, $(L _{m-j}, v _{j})$ is a
Henselian discrete valued field with $\widehat L _{m-j} = L
_{m-j-1}$, and by Lemma \ref{lemm4.4} and
\cite[Corollary~8.2.40]{Li}, $O _{v _{j}}(L _{m-j})$ is an
excellent ring. Since $K _{0}$ does not admit proper algebraic
extensions in $K _{m}$, this implies $p$ does not divide the
degree of any finite extension of $K _{m}$ in $L _{m}$. Therefore,
it follows from Lemmas \ref{lemm3.1} and \ref{lemm3.2} that
abrd$_{p}(L _{0}(X)) = {\rm abrd}_{p}(K _{0}(X)) < \infty $, and
abrd$_{p'}(L _{0}(X)) = 0$, for every prime $p' \neq p$.
Similarly, it can be deduced from Lemma \ref{lemm3.2} that
abrd$_{p}(L _{m}) = {\rm abrd}_{p}(K _{m})$. Thus the inequalities
abrd$_{p}(L _{m}(X)) < \infty $ and $b _{p}(K _{m}) \le b _{p}(K
_{0}) + m$ reduce to consequences of \cite[Corollaries~5.7,
5.8]{HHKr} and \cite[Corollary~1.3]{Li11}. As $p$ is an arbitrary
element of $\mathbb{P}_{q}$ and $K _{m}$ is virtually perfect, the
obtained result and Lemma \ref{lemm2.1} show that $K = K _{m}(X)$
is a $\Phi _{\rm Br}$-field, so Theorem \ref{theo2.2} is proved.
\par
\medskip
Note finally that, for any $m \in \mathbb{N}$, Theorem
\ref{theo2.2} leaves open the question of whether a field is
Brauer finite-dimensional, provided that it is a
finitely-generated extension of an $m$-local field $K _{m}$ with
char$(K _{m}) = 0$ and an $m$-th residue field $K _{0}$ such that
char$(K _{0}) = q > 0$ and ddim$(K _{0})$ is finite. The answer
will be positive if it turns out that ddim$(K _{m,p}) < \infty
\colon p \in \mathbb{P}$, where $K _{m,p}$ is the fixed field of
some Sylow pro-$p$ subgroup of $\mathcal{G}_{K_{m}}$, for each
$p$.

\end{document}